\documentclass[a4paper]{article}
\usepackage{theorem,amsmath,amssymb,amscd}
\usepackage[small]{titlesec}

\theoremstyle{change} \allowdisplaybreaks 

\newcommand{\R}{\mathbb{R}}
\newcommand{\C}{\mathbb{C}}
\newcommand{\Z}{\mathbb{Z}}
\newcommand{\N}{\mathbb{N}}
\newcommand{\Q}{\mathbb{Q}}
\newcommand{\A}{\mathbb{A}}

\newcommand{\Kl}{\mathrm{Kl}}
\newcommand{\St}{\mathrm{St}}
\newcommand{\GL}{\mathrm{GL}}
\newcommand{\GU}{\mathrm{GU}}

\newcommand{\GSp}{\mathrm{GSp}}
\newcommand{\SSp}{\mathrm{Sp}}
\newcommand{\OF}{\mathfrak{o}}
\newcommand{\p}{\mathfrak{p}}
\renewcommand{\P}{\mathfrak{P}}

\newcommand{\tr}{{\rm tr}}
\newcommand{\AI}{\mathcal{AI}}
\newcommand{\mat}[4]{{\setlength{\arraycolsep}{0.5mm}\left[
\begin{array}{cc}#1&#2\\#3&#4\end{array}\right]}}
\newcommand{\qed}{\hspace*{\fill}\rule{1ex}{1ex}}
\def\qdots{\mathinner{\mkern1mu\raise0pt\vbox{\kern7pt\hbox{.}}\mkern2mu
\raise3.4pt\hbox{.}\mkern2mu\raise7pt\hbox{.}\mkern1mu}}

\newcommand{\nl}{

\vspace{2ex}}
\newcommand{\nll}{

\vspace{1ex}}

\newtheorem{thm}{Theorem.}[section]
\newtheorem{theorem}{Theorem.}[section]

\newtheorem{lemma}[thm]{Lemma.}

\newtheorem{proposition}[thm]{Proposition.}

\begin{document}

\begin{center}
{\Large $L$-functions for $\GSp_4\times \GL_2$ in the case of high $\GL_2$ conductor}

\vspace{2ex} Ameya Pitale\footnote{Department of Mathematics, University of Oklahoma,
Norman, OK 73019-0315, {\tt ameya@math.ou.edu}},
Ralf Schmidt\footnote{Department of Mathematics, University of Oklahoma,
Norman, OK 73019-0315, {\tt rschmidt@math.ou.edu}}

\vspace{4ex}
\begin{minipage}{80ex}
 {\small {\sc Abstract.} Furusawa \cite{Fu} has given an integral representation for
 the degree $8$ L-function of $\GSp_4\times\GL_2$ and has carried out the
 unramified calculation. The local $p$-adic zeta integrals were calculated
 in the work \cite{P-S} under the assumption that the $\GSp_4$ representation $\pi$
 is unramified and the $\GL_2$ representation $\tau$ has conductor $\p$. In the
 present work we generalize to the case where the $\GL_2$
 representation has arbitrarily high conductor. The result is that the
 zeta integral represents the local Euler factor $L(s,\pi\times\tau)$ in all cases.
 As a global application we obtain a special value result for
 a $\GSp_4\times\GL_2$ global $L$-function coming from classical holomorphic
 cusp forms with arbitrarily high level for the elliptic modular form.}
\end{minipage}
\end{center}
\section{Introduction}
Let $\pi = \otimes \pi_{\nu}$ and $\tau = \otimes \tau_{\nu}$ be
irreducible, cuspidal, automorphic representations of $\GSp_4(\A)$
and $\GL_2(\A)$, respectively. Here, $\A$ is the ring of adeles of
a number field $F$. We investigate the degree eight
twisted $L$-functions $L(s, \pi \times \tau)$ of $\pi$ and $\tau$,
which are important for a number of reasons. For
example, when $\pi$ and $\tau$ are obtained from holomorphic
modular forms, then Deligne \cite{De1} has conjectured that a
finite set of special values of $L(s, \pi \times \tau)$ are
algebraic up to certain period integrals. Another very important
application is the conjectured Langlands functorial transfer of
$\pi$ to an automorphic representation of $\GL_4(\A)$. One
approach to obtain the transfer to $\GL_4(\A)$ is to use the
converse theorem due to Cogdell and Piatetski-Shapiro \cite{CPS},
which requires precise information about the $L$-functions $L(s, \pi
\times \tau)$.
 
In \cite{Fu}, Furusawa has obtained an integral representation for $L(s, \pi \times \tau)$ in the special case where $\pi$ and $\tau$ correspond to holomorphic cusp forms of full level and same weight. In \cite{P-S}, we extended Furusawa's integral representation to the case where $\tau$ corresponds to a cusp form (holomorphic or non-holomorphic) of square-free level. At the non-archimedean place $\nu$ dividing the level, this means that $\tau_\nu$ is mildly ramified, namely, it is an unramified twist of the Steinberg representation of $\GL_2$.

In this paper, we will compute the local non-archimedean integral for any irreducible, admissible, generic representation $\tau_\nu$ with unramified central character; the
conductor of $\tau_\nu$ can be arbitrarily high. In the local case, we obtain the following result.

{\bf Theorem 1.}
 \emph{Let $F_{\nu}$ be a non-archimedean local field with characteristic zero. Let
 $\pi_{\nu}$ be an unramified, irreducible, admissible
 representation of $\GSp_4(F_{\nu})$. Let $\tau_{\nu}$ be an irreducible, admissible, generic representation of $\GL_2(F_\nu)$ with unramified central character and conductor $\p^n, n \geq 2$. Then we can make a choice of
 vectors such that the local integral (defined in (\ref{local-zeta-integral}))
 is given by
 $$
  Z_{\nu}(s) = \frac{q-1}{q^{3(n-1)}(q+1)(q^4-1)} (1-\Big(\frac L{\p}\Big)q^{-1})q^n.
 $$
 Here, $q$ is the cardinality of the residue class field of $F_{\nu}$ and $L_{\nu}$ is the degree $2$ extension of $F_\nu$ defined in the next section.
}

The $n=0$ case was done in \cite{Fu} and the $n=1$ case was done in \cite{P-S}. Note that, for representations $\pi_\nu$ and $\tau_\nu$ as described in the above theorem, we have $L(s, \pi_\nu \times \tau_\nu) = 1$, and hence the integral $Z_{\nu}(s)$ indeed computes the $L$-function up to a constant. Not surprisingly, the case $n > 1$ is not a straightforward generalization of the case $n=1$, but requires different arguments.
Making the ``correct'' choice of local vectors to be used to compute the local integral is delicate and, probably, is the main contribution of this paper. For example, we will
have to make a choice of local compact subgroup $K^\#(\P^n)$, for which the Borel congruence
subgroup turns out to be too small, while the Klingen congruence subgroup is too large;
the group we will work with lies in between these two natural congruence subgroups.
We would like to point out that, in general, the ramified calculation is very complicated in most situations and is rarely carried out in much of the available work on integral representations of $L$-functions.

We will now describe the global case. 

{\bf Theorem 2.}
\emph{
Let $\Phi$ be a Siegel cuspidal eigenform of weight $l$ with respect to $\SSp_4(\Z)$ satisfying the two assumptions described in Section \ref{global-section}. Let $N$ be any positive integer. Let $f$ be a Maa{\ss}
 Hecke eigenform of weight $l_1 \in \Z$ with respect to $\Gamma_0(N)$. If $f$ lies in a
 holomorphic discrete series with lowest weight $l_2$, then assume that $l_2 \leq l$. Let $\pi_{\Phi}$ and $\tau_f$ be the corresponding cuspidal
 automorphic representations of $\GSp_4(\A_{\Q})$ and
 $\GL_2(\A_{\Q})$, respectively. Then a choice of local vectors can be made
 such that the global integral $Z(s)$ defined in (\ref{global integral calculation}) is given by
$$
  Z(s) = \kappa_{\infty} \kappa_N
  \frac{L(3s+\frac 12, \pi_{\Phi} \times \tau_f)}{\zeta(6s+1)
  L(3s+1, \tau_f\times\AI(\Lambda))},
 $$
 where $\kappa_{\infty}, \kappa_N$ are obtained from the local computations, $\Lambda$ is the Bessel character defined in the next section and $\AI(\Lambda)$ is the representation of $\GL_2(\A)$ obtained by automorphic induction from the character $\Lambda$.
 }
 
We can obtain a special value result in the spirit of Deligne's conjectures from the above theorem. For this we need to make an additional assumption on $N$. Suppose $D$ is the positive integer as in the Assumption 1 on $\Phi$ and $p$ is a prime that divides $D$, then we assume that the highest power of $p$ dividing $N$ is not equal to $2$.

{\bf Theorem 3.}
\emph{ Let $\Phi$ be a Siegel cuspidal eigenform of weight $l$ with respect to $\SSp_4(\Z)$. Let $N$ be any positive integer satisfying the condition described above. Let $\Psi$ be a holomorphic cuspidal eigenform of weight $l$ with respect to $\Gamma_0(N)$. Then 
$$ \frac{L(\frac l2 - 1, \pi_{\Phi} \times \tau_{\Psi})}{\pi^{5l-8}
  \langle \Phi,\Phi \rangle \langle \Psi,\Psi \rangle_1} \in \bar{\Q}.$$ }

We refer the reader to \cite{BH}, \cite{Fu}, \cite{P-S} and \cite{Sa} for related results on special values of the $\GSp_4 \times \GL_2$ $L$-functions. 

This paper is organized as follows. In Section \ref{preliminaries-section}, we describe the notations and general setup for the non-archimedean calculation. This is just a brief outline of a more detailed description given in \cite{P-S}. In Section \ref{compact-grp-sec}, we describe the choice of the local vector used for the integral calculation. The computation of the integral basically involves three main things -- a suitable double coset decomposition, computation of the volumes of the double cosets and computing the expression obtained by substituting the local vectors in the integral. This is carried out in Sections \ref{double coset decomposition sec}, \ref{volume-sec} and \ref{main-thm-sec}. In Section \ref{global-section}, we put together the local non-archimedean result with the computations from \cite{Fu} and \cite{P-S} to obtain the global theorem and the special value result.

\section{Preliminaries}\label{preliminaries-section}
In this section, we will briefly recall the setup and terminology in the non-archimedean local case from \cite{P-S}.
\begin{enumerate}
\item Let $F$ be a non-archimedean local field of characteristic zero.
Let $\OF$, $\p$, $\varpi$, $q$ be the ring of integers, prime ideal, uniformizer
and cardinality of the residue class field $\OF/\p$, respectively.

\item Fix three
elements $a,b,c\in F$ such that $d := b^2-4ac \neq 0$. Let $
 L = F(\sqrt{d})$ if $d \notin F^{\times2}$ and $L = F \oplus F$ if $d \in F^{\times2}$. If $L$ is a field,
 then let $\OF_L$ be its ring of integers. If $L = F \oplus F$,
then let $\OF_L = \OF \oplus \OF$. If $L$ is a field, we denote by
$\bar x$ the Galois conjugate of $x\in L$ over $F$. If $L=F\oplus
F$, let $\overline{(x,y)}=(y,x)$.

\item We shall assume that $a,b\in\OF$ and $c\in\OF^\times$. In
addition, we will assume that if $d \not\in F^{\times2}$, then $d$
is the generator of the discriminant of $L/F$ and
  if $d \in F^{\times2}$, then $d \in \OF^{\times}$.

\item Let
$
 \alpha:=\frac{b+\sqrt{d}}{2c}$ if $L$ is a field and  $ \alpha :=
 \Big(\frac{b+\sqrt{d}}{2c},\frac{b-\sqrt{d}}{2c}\Big)$ if
$L=F\oplus F$. Let $\eta
=\begin{bmatrix}1&0&&\\\alpha&1&&\\&&1&-\bar{\alpha}\\&&0&1\end{bmatrix}$.

\item We fix the following ideal in $\OF_L$,
\begin{equation}\label{ideal defn}\renewcommand{\arraystretch}{1.3}
 \P := \p\OF_L = \left\{
                  \begin{array}{l@{\qquad\text{if }}l}
                    \p_L & \big(\frac L{\p}\big) = -1,\\
                    \p_L^2 & \big(\frac L{\p}\big) = 0,\\
                    \p \oplus \p & \big(\frac L{\p}\big) = 1.
                  \end{array}
                \right.
\end{equation}
Here, $\big(\frac L{\p}\big)$ is the Legendre symbol defined in
(22) of \cite{P-S} and $\p_L$ is the maximal ideal of $\OF_L$ when
$L$ is a field extension. Note that $\P$ is prime only if
$\big(\frac L\p\big)=-1$. We have $\P^n\cap\OF=\p^n$ for all
$n\geq0$.

\item Let \begin{align*}
 H(F) = \GSp_4(F) &:= \{g \in \GL_4(F) : \, ^tgJg = \mu(g)J,\:\mu(g)\in F^{\times} \}, \\
 G(F) = \GU(2,2;L) &:= \{g \in \GL_4(L) : \, ^t\bar{g}Jg = \mu(g)J,\:\mu(g)\in F^{\times}\},
\end{align*}
where $J = \mat{}{1_2}{-1_2}{}$. Let the subgroups $P, M^{(1)},
M^{(2)}$ and $N$ of $G(F)$ be as defined in Sect.\ 2.1 of
\cite{P-S}.

\item For $S = \mat{a}{b/2}{b/2}{c}$ define
$T(F)=\{g\in\GL_2(F):\:^tgSg=\det(g)S\} = \{
\mat{x+yb/2}{yc}{-ya}{x-yb/2} : x,y \in F\} \simeq L^\times$. Let
$U(F)=\{\mat{1_2}{X}{}{1_2}\in\GSp_4(F):\:^tX=X\}$. Then the
Bessel subgroup $R(F)$ is defined as $R(F) = T(F) U(F)$.

 \item Let $K^H=\GSp_4(\OF)$.
  From (3.4.2) of \cite{Fu}, we have the disjoint double coset decomposition
  \begin{equation}\label{RFKHrepresentativeseq}
   H(F)=\bigsqcup_{l\in\Z}\bigsqcup_{m\geq0}R(F)h(l,m)K^H, \quad\mbox{ where }
   h(l,m)=\begin{bmatrix}\varpi^{2m+l}\\&\varpi^{m+l}\\&&1\\&&&\varpi^m\end{bmatrix}.
  \end{equation}

\item Fix a nontrivial character $\psi$ of $F$ with conductor
$\OF$. Let $\theta$ be a character of $U(F)$ obtained from $\psi$
and $S$ as in Sect.\ 2.2 of \cite{P-S}. Let $\Lambda$ be any
character of $T(F) \simeq L^\times$. Then we get a character
$\Lambda \otimes \theta$ of $R(F)$ as in Sect.\ 2.2 of \cite{P-S}.
Let $(\pi, V_{\pi})$ be an unramified, irreducible, admissible
representation of $H(F)$ with central character $\omega_\pi$.
Assume that $\Lambda|_{F^\times} = \omega_\pi$. We assume that
$V_\pi$ is a Bessel model for $\pi$ with respect to the character
$\Lambda \otimes \theta$ of $R(F)$. Let $B$ denote the unique
spherical vector in $V_\pi$ satisfying $B(1) = 1$. By Lemma 3.4.4
of \cite{Fu} we have $B(h(l,m))=0$ for $l<0$.

 \item Let $(\tau, V_\tau)$ be any generic, irreducible, admissible
  representation of $\GL_2(F)$ with an unramified central character
  $\omega_\tau$. Let $\chi_0$ be a character of $L^\times$ such that
  $\chi_0|_{F^\times} = \omega_\tau$. Let $\chi$ be another
  character of $L^\times$ such that $\chi(\zeta) =
  \Lambda(\bar{\zeta})^{-1} \chi_0(\bar{\zeta})^{-1}$. Let
  $I(s,\chi, \chi_0, \tau)$ be the induced representation of $G(F)$
  constructed in Sect.\ 2.3 of \cite{P-S}. Further below we will
  construct a function $W^\#\in I(s,\chi, \chi_0, \tau)$. Our main
  local result will be the evaluation of the integral
  \begin{equation}\label{local-zeta-integral}
   Z(s) = \int\limits_{R(F) \backslash H(F)} W^\#(\eta h, s) B(h) dh.
  \end{equation}
\end{enumerate}
\section{The local compact group and the function $W^\#$}\label{compact-grp-sec}
We define congruence subgroups of $\GL_2(F)$, as follows. For
$n=0$ let $K^{(1)}(\p^0)=\GL_2(\OF)$. For $n>0$ let
\begin{equation}\label{K1defeq}
 K^{(1)}(\p^n)=\GL_2(F)\cap\mat{\OF^\times}{\OF}{\p^n}{\OF^\times}.
\end{equation}
Note that we are only considering representations $\tau$ with
unramified central character. For such representations, the conductor
(or level) of $\tau$ can be defined in terms of the congruence
subgroups (\ref{K1defeq}). More precisely,
if $n$ is the smallest integer for which $V_\tau$ has a vector
that is invariant under $K^{(1)}(\p^n)$, then we say that
$\p^n$ is the conductor of $\tau$. The space of such
invariant vectors is one dimensional. We assume that $V_\tau$ is
the Whittaker model of $\tau$ with respect to the character of $F$
given by $\psi^{-c}(x)=\psi(-cx)$. Note that this character has
conductor $\OF$ by our assumptions on $\psi$ and $c$. Let $W^{(0)}\in V_\tau$ be the local
newform, i.e., the essentially unique non-zero $K^{(1)}(\p^n)$
invariant vector in $V_\tau$. We can make it unique by requiring
that $W^{(0)}(1)=1$.

Let $\P = \p \OF_L$, as above. Let
\begin{equation}\label{iwahori defn}
 I := \{g \in \GU(2,2;\OF_L):\:g\equiv\begin{bmatrix}\ast&0&\ast&\ast\\
 \ast&\ast&\ast&\ast\\0&0&\ast&\ast\\ 0&0&0&\ast\end{bmatrix}\pmod{\P} \}
\end{equation}
be the Iwahori subgroup and
\begin{equation}\label{Klingen-congruence defn}
 \Kl(\P^n) := \{g \in \GU(2,2;\OF_L):\:g \equiv \begin{bmatrix}\ast&0&\ast&\ast \\
  \ast&\ast&\ast&\ast\\\ast&0&\ast&\ast\\ 0&0&0&\ast\end{bmatrix} \pmod{\P^n} \}
\end{equation}
be the Klingen congruence subgroup. We define $K^\#(\P^0) :=
\GU(2,2;\OF_L)$ and for $n \geq 1$
\begin{equation}\label{U-compact defn}
 K^\#(\P^n) := I \cap \Kl(\P^n) = \GU(2,2;\OF_L) \cap
  \begin{bmatrix}\OF_L^{\times} & \P^n & \OF_L & \OF_L \\
  \OF_L & \OF_L^{\times} & \OF_L & \OF_L \\ \P & \P^n & \OF_L^{\times} & \OF_L \\
  \P^n & \P^n & \P^n & \OF_L^{\times} \end{bmatrix}.
\end{equation}
Furthermore, let
\begin{equation}\label{U-compact defn2}
 K^\#(\p^n) := K^\#(\P^n) \cap \GSp_4(F) = \GU(2,2;\OF) \cap
  \begin{bmatrix}\OF^{\times} & \p^n & \OF & \OF \\
  \OF & \OF^{\times} & \OF & \OF \\ \p & \p^n & \OF^{\times} & \OF \\
  \p^n & \p^n & \p^n & \OF^{\times} \end{bmatrix}.
\end{equation}
Note that $K^\#(\P)=I$. Also, note that $K^\#(\P^n)$ is slightly different
from the group (with the same name) defined in Sect.\ 3.3 of
\cite{P-S}.

We will now define the specific function $W^\#$ which we will use
to evaluate the integral (\ref{local-zeta-integral}). Since a similar
definition has been made in Sect.\ 2.3 of \cite{P-S}, we will omit
some details. We first extend $W^{(0)}$ to a function on $M^{(2)}(F)$ via
$W^{(0)}(ag)=\chi_0(a)W^{(0)}(g)$ for $a\in L^\times$ and $g\in\GL_2(F)$.
Given a complex number $s$, there exists a unique function
$W^\#(\,\cdot\,,s):\:G(F)\rightarrow\C$ with the following properties.
\begin{enumerate}
 \item If $g\notin M(F)N(F)K^\#(\P^n)$, then $W^\#(g,s)=0$.
 \item If $g=mnk$ with $m\in M(F)$, $n\in N(F)$, $k\in K^\#(\P^n)$, then
  $W^\#(g,s)=W^\#(m,s)$.
 \item For $\zeta\in L^\times$ and $\mat{\alpha}{\beta}{\gamma}{\delta}\in M^{(2)}(F)$,
  \begin{equation}\label{Wsharpformulaeq}
   W^\#(\begin{bmatrix}\zeta\\&1\\&&\bar{\zeta}^{-1}\\&&&1\end{bmatrix}
   \begin{bmatrix}1\\&\alpha&&\beta\\&&\mu\\&\gamma&&\delta\end{bmatrix},s)
   =|N(\zeta)\cdot\mu^{-1}|^{3(s+1/2)}\chi(\zeta)\,
   W^{(0)}(\mat{\alpha}{\beta}{\gamma}{\delta}).
  \end{equation}
  Here $\mu=\bar\alpha\delta-\beta\bar\gamma$.
\end{enumerate}

As described in Sect.\ 3.5 of \cite{P-S}, the integral
(\ref{local-zeta-integral}) reduces to
\begin{equation}\label{local-zeta-integral-2}
 Z(s)=\sum_{l,m\geq0}\;\sum_i\;B(h(l,m))\,W^\#(\eta h(l,m)s_i,s)
  \int\limits_{K_{l,m}\backslash K_{l,m}s_iK^\#(\p^n)}\,dh,
\end{equation}
where $K_{l,m}:=h(l,m)^{-1}R(F)h(l,m)\cap K^H$ and $\{s_i\}$ is a
system of representatives for the double coset space
$K_{l,m}\backslash K^H/K^\#(\p^n)$. We will now follow the three
steps outlined in Sect.\ 3.5 of \cite{P-S} to obtain a suitable
subset $\{s_i'''\}$ of $\{s_i\}$ for which $W^\#(\eta h(l,m)s_i,s)\neq0$.
\section{Double coset decomposition}\label{double coset decomposition sec}
\subsection{The cosets $K^\#(\p^0)/K^\#(\p^n)$}\label{KsharpcosetssecU}
We need to determine representatives for the coset space
\begin{equation}\label{newK2KHcosetseq12U}
 K^\#(\p^0)/K^\#(\p^n),\qquad\text{where }K^\#(\p^0)=K^H=\GSp(4,\OF).
\end{equation}
Note that this coset space is isomorphic to
\begin{equation}\label{newK2KHcosetseq22U}
 K^\#_1(\p^0)/K^\#_1(\p^n),\qquad\text{where }
 K^\#_1(\p^n)=K^\#(\p^n)\cap\{g\in H(F):\mu(g)=1\}.
\end{equation}
Let
\begin{equation}\label{Weylgroupelementsdefeq3U}
 \qquad s_1=\begin{bmatrix}&1\\1\\&&&1\\&&1\end{bmatrix},
 \qquad s_2=\begin{bmatrix}&&1\\&1\\-1\\&&&1\end{bmatrix}.
\end{equation}
It follows from the Bruhat decomposition for $\SSp(4,\OF/\p)$ that
\begin{align}
 \label{newK2KHcosetseq41b}K^\#(\p^0)&=K^\#(\p^1) \sqcup\:\bigsqcup_{x\in\OF/\p}
  \begin{bmatrix}1\\x&1\\&&1&-x\\&&&1\end{bmatrix}s_1K^\#(\p^1) \sqcup\:\bigsqcup_{x\in\OF/\p}
  \begin{bmatrix}1&&x\\&1\\&&1&\\&&&1\end{bmatrix}s_2K^\#(\p^1)\\
 \label{newK2KHcosetseq44b}&\sqcup\:\bigsqcup_{x,y\in\OF/\p}
  \begin{bmatrix}1\\x&1&&y\\&&1&-x\\&&&1\end{bmatrix}s_1s_2K^\#(\p^1) \sqcup\:\bigsqcup_{x,y\in\OF/\p}
  \begin{bmatrix}1&&x&y\\&1&y\\&&1\\&&&1\end{bmatrix}s_2s_1K^\#(\p^1)\\
 \label{newK2KHcosetseq46b}&\sqcup\:\bigsqcup_{x,y,z\in\OF/\p}
  \begin{bmatrix}1&&&y\\x&1&y&xy+z\\&&1&-x\\&&&1\end{bmatrix}s_1s_2s_1K^\#(\p^1) \sqcup\:\bigsqcup_{x,y,z\in\OF/\p}
  \begin{bmatrix}1&&x&y\\&1&y&z\\&&1\\&&&1\end{bmatrix}s_2s_1s_2K^\#(\p^1)\\
 \label{newK2KHcosetseq48b}&\sqcup\:\bigsqcup_{w,x,y,z\in\OF/\p}
  \begin{bmatrix}1&&x&y\\w&1&wx+y&wy+z\\&&1&-w\\&&&1\end{bmatrix}s_1s_2s_1s_2K^\#(\p^1).
\end{align}
Let $n\geq1$. It is easy to see that
\begin{equation}\label{kgeq2Ksharpdecompeq2U}
 K^\#(\p^1)=\bigsqcup\limits_{w,y,z\in\OF/\p^{n-1}}
 \begin{bmatrix}1&w\varpi\\&1\\&&1\\&&-w\varpi&1\end{bmatrix}
  \begin{bmatrix}1\\&1\\&y\varpi&1\\y\varpi&z\varpi&&1\end{bmatrix}K^\#(\p^n).
\end{equation}
Let $\{r_i\}$ be the system of representatives for $K^\#(\p^0)/K^\#(\p^1)$ determined
in (\ref{newK2KHcosetseq41b}) -- (\ref{newK2KHcosetseq48b}). Combining these with
(\ref{kgeq2Ksharpdecompeq2U}) we get
\begin{equation}\label{kgeq2KsharpKHdecompeq2U}
 K^H=\bigsqcup\limits_i \bigsqcup\limits_{w,y,z\in\OF/\p^{n-1}}r_i \begin{bmatrix}1&w\varpi\\&1\\&&1\\&&-w\varpi&1\end{bmatrix}
  \begin{bmatrix}1\\&1\\&y\varpi&1\\y\varpi&z\varpi&&1\end{bmatrix}
 K^\#(\p^n).
\end{equation}
Recall that we are interested in the double cosets $K_{l,m}\backslash K^H/K^\#(\p^n)$, where
$K_{l,m}=h(l,m)^{-1}R(F)h(l,m)\cap K^H$.
\subsection{Step $1$: Preliminary decomposition}
Observe that $K_{l,m}$ contains all elements
$\begin{bmatrix}1&&\OF&\OF\\&1&\OF&\OF\\&&1\\&&&1\end{bmatrix}$.
From (\ref{kgeq2KsharpKHdecompeq2U}) we get the following preliminary
decomposition, which is not disjoint.
\begin{eqnarray}
 K^H&=&  \bigcup\limits_{y,z,w\in\OF/\p^{n-1}}K_{l,m}
   \begin{bmatrix}1&w\varpi\\&1\\&&1\\&&-w\varpi&1\end{bmatrix}
   \begin{bmatrix}1\\&1\\&y\varpi&1\\y\varpi&z\varpi&&1\end{bmatrix}
   K^\#(\p^n) \label{prelimnewK2KHcosetseq41U}\\
 &&\cup\bigcup\limits_{w\in\OF/\p^n}\;
   \bigcup\limits_{y,z\in\OF/\p^{n-1}}K_{l,m}
   \begin{bmatrix}1\\w&1\\&&1&-w\\&&&1\end{bmatrix}
   \begin{bmatrix}1\\&1\\z\varpi&y\varpi&1\\y\varpi&&&1\end{bmatrix}
   s_1K^\#(\p^n) \label{prelimnewK2KHcosetseq42U}\\
 &&\cup\bigcup\limits_{w,y,z\in\OF/\p^{n-1}}\;K_{l,m}
   \begin{bmatrix}1\\&1\\&w\varpi&1&\\w\varpi&&&1\end{bmatrix}
   \begin{bmatrix}1&y\varpi\\&1\\&&1\\&z\varpi&-y\varpi&1\end{bmatrix}
   s_2K^\#(\p^n) \label{prelimnewK2KHcosetseq43U}\\
 &&\cup\bigcup\limits_{w\in\OF/\p^n}\;\bigcup\limits_{y,z\in\OF/\p^{n-1}}
   \;K_{l,m}\begin{bmatrix}1\\w&1\\&&1&-w\\&&&1\end{bmatrix}
   \begin{bmatrix}1&\\&1\\y\varpi&z\varpi&1\\z\varpi&&&1\end{bmatrix}
   s_1s_2K^\#(\p^n) \label{prelimnewK2KHcosetseq44U}\\
 &&\cup\bigcup\limits_{w\in\OF/\p^{n-1}}K_{l,m}
   \begin{bmatrix}1&w\varpi\\&1\\&&1\\&&-w\varpi&1\end{bmatrix}
   s_2s_1K^\#(\p^n) \label{prelimnewK2KHcosetseq45U}\\
 &&\cup\bigcup\limits_{w\in\OF/\p^n}K_{l,m}
   \begin{bmatrix}1\\w&1\\&&1&-w\\&&&1\end{bmatrix}s_1s_2s_1K^\#(\p^n)
   \label{prelimnewK2KHcosetseq46U}\\
 &&\cup\bigcup\limits_{w\in\OF/\p^{n-1}}K_{l,m}
   \begin{bmatrix}1&w\varpi\\&1\\&&1&\\&&-w\varpi&1\end{bmatrix}s_2s_1s_2K^\#(\p^n)
   \label{prelimnewK2KHcosetseq47U}\\
 &&\cup\bigcup\limits_{w\in\OF/\p^n}K_{l,m}
   \begin{bmatrix}1\\w&1\\&&1&-w\\&&&1\end{bmatrix}s_1s_2s_1s_2K^\#(\p^n).
   \label{prelimnewK2KHcosetseq48U}
\end{eqnarray}
\subsection{Step $2$: Support of $W^\#$}\label{support of W-U}
We assumed that $c\in\OF^\times$, so that $\alpha\in\OF_L$.
We have $\eta h(l,m) = h(l,m)\eta_m$, where for $m\geq0$ we let
\begin{equation}\label{etamdefeq2U}
 \eta_m=\begin{bmatrix}1\\\alpha\varpi^m&1\\&&1&-\bar\alpha\varpi^m\\&&&1\end{bmatrix}.
\end{equation}
Fix $l,m\geq0$, and let $r$ run through the representatives for
$K_{l,m}\backslash K^H/K^\#(\p^n)$ from
(\ref{prelimnewK2KHcosetseq41U}) -- (\ref{prelimnewK2KHcosetseq48U}).
In view of (\ref{local-zeta-integral-2})
we want to find out for which $r$ is $\eta h(l,m)r \in M(F)N(F)
K^\#(\P^n)$, since this set is the support of $W^\#$. Since
$h(l,m)\in M(F)$, this is equivalent to $\eta_m r \in M(F)N(F)
K^\#(\P^n)$. Hence, this condition depends only on $m\geq0$ and
not on the integer $l$.
\begin{enumerate}
 \item Let
  $
   r= \begin{bmatrix}1&w\varpi\\&1\\&&1\\&&-w\varpi&1\end{bmatrix}  \begin{bmatrix}1\\&1\\&y\varpi&1\\y\varpi&z\varpi&&1\end{bmatrix}
  $
  with $w,y,z\in\OF/\p^{n-1}$.
Suppose $\eta_m r = \tilde{m}\tilde{n}k$ with $\tilde{m} \in M(F)$, $\tilde{n} \in N(F)$ and $k \in K^\#(\P^n)$. Let $A = (\tilde{m}\tilde{n})^{-1}\eta_m r$. Looking at the $(3,2)$ and $(3,3)$ coefficient of $A$ we get
$$
 y+\varpi^{m+1} \alpha w y - \varpi^m \alpha z \in \P^{n-1},\qquad\text{and hence}
 \qquad\alpha \varpi^m (\varpi w y - z) + y \in \P^{n-1}.
$$
If $\nu(\varpi w y - z) < n-m-1$ then $\alpha + y/(\varpi^m(\varpi
w y - z)) \in \P$, which contradicts Lemma 3.1.1 (ii) of \cite{P-S}.
Hence, $\nu(\varpi w y - z) \geq n-m-1$, which implies
$\varpi^m(\varpi w y - z) \in \p^{n-1}$. It follows that $y \in\p^{n-1}$.
To summarize, necessary conditions for $A\in K^\#(\P^n)$
are $y=0$ and $z \in (\p^{n-m-1} \cap\OF)/\p^{n-1}$.
The following matrix identity shows that these are also sufficient conditions.
\begin{eqnarray}\label{type-1-relevance-eqn}
 \eta_m r &=& \begin{bmatrix}a^{-1}&&&\\&a&&\\&&\bar{a}&\\
 &\varpi z \bar{a}^{-1}&&\bar{a}^{-1}\end{bmatrix} \begin{bmatrix}1&\varpi w a&&\\
 &1&&\\&&1&\\&&-\varpi w \bar{a}&1\end{bmatrix} \nonumber\\
 && \qquad \begin{bmatrix}1&&&\\ \varpi^m\alpha a^{-1}&1&&\\
 &-\varpi^{m+1}\bar{\alpha}z\bar{a}^{-1}&1&-\varpi^m\bar{\alpha} \bar{a}^{-1}\\
 -\varpi^{m+1}\alpha za^{-1}&&&1\end{bmatrix} \in M(F)N(F)K^\#(\P^n),
\end{eqnarray}
where $a = 1 + \varpi^{m+1} \alpha w \in \OF_L^{\times}$. Hence, the values of $w,y,z$ for which $\eta_m r \in M(F)N(F)K^\#(\P^n)$ are
  $$
   w \in \OF/\p^{n-1},\; y=0,\; z \in (\p^{n-m-1} \cap \OF)/\p^{n-1}.
  $$
\item Let $r=  \begin{bmatrix}1\\w&1\\&&1&-w\\&&&1\end{bmatrix}
  \begin{bmatrix}1\\&1\\z\varpi&y\varpi&1\\y\varpi&&&1\end{bmatrix}
  s_1$ with $w\in\OF/\p^n$ and $y,z\in\OF/\p^{n-1}$.
 Suppose $\eta_m r = \tilde{m}\tilde{n}k$ with $\tilde{m} \in M(F), \tilde{n} \in N(F)$
 and $k \in K^\#(\P^n)$. Let $A = (\tilde{m}\tilde{n})^{-1}\eta_m r$.
 Looking at the $(3,2)$ and $(3,3)$ coefficients of $A$ we get
 $$
  \beta:=\varpi^m \alpha + w \in \OF_L^{\times}\qquad\mbox{and}\qquad
  \varpi^m \alpha y + wy - z\in \P^{n-1}.
 $$
 If $\nu(y) < n-m-1$, then $\alpha + (wy - z)/(\varpi^m y) \in \P$,
 which contradicts Lemma 3.1.1(ii) of \cite{P-S}.
 Hence, $\nu(y) \geq n-m-1$, which implies $wy - z \in \P^{n-1}$.
 We may therefore assume that $z = wy$. To summarize, necessary
 conditions for $A \in K^\#(\P^n)$ are $\varpi^m \alpha + w \in
 \OF_L^{\times}$, $y \in (\p^{n-m-1} \cap \OF)/\p^{n-1}$ and
 $z = wy$. The following matrix identity shows that these are also sufficient conditions.
 \begin{eqnarray}
 \eta_m r &=& \begin{bmatrix}-\beta^{-1}&&&\\&\beta&&\\&&-\bar{\beta}&\\&\varpi w y \bar{\beta}^{-1}&&\bar{\beta}^{-1}\end{bmatrix} \begin{bmatrix}1&-\beta&&\\&1&&\\&&1&\\&&\bar{\beta}&1\end{bmatrix} \nonumber \\
 && \qquad \begin{bmatrix}1&&&\\ \beta^{-1}&1&&\\-\varpi y \bar{\beta}^{-1}& \varpi^{m+1} \bar{\alpha} y \bar{\beta}^{-1}&1&-\bar{\beta}^{-1}\\ \varpi^{m+1} \alpha y \beta^{-1} &&&1\end{bmatrix} \in M(F)N(F)
K^\#(\P^n). \label{type-2-relevance-eqn}
\end{eqnarray}
 Hence, the values of $w,y,z$ for which $\eta_m r \in M(F)N(F)K^\#(\P^n)$ are as follows.
  \begin{enumerate}
   \item If $m=0$, then all $w \in \OF/\p^n$ such that $\alpha + w \in \OF_L^{\times}$
    and $y=z=0$.
   \item If $m>0$, then all $w \in \OF^{\times}$, $y \in (\p^{n-m-1} \cap \OF)/\p^{n-1}$
    and $z=wy$.
  \end{enumerate}

 \item Let $r=\begin{bmatrix}1\\&1\\&w\varpi&1&\\w\varpi&&&1\end{bmatrix}
  \begin{bmatrix}1&y\varpi\\&1\\&&1\\&z\varpi&-y\varpi&1\end{bmatrix}s_2$
  with $w,y,z\in\OF/\p^{n-1}$.

  Then $\eta_mr\notin M(F)N(F)K^\#(\P^n)$,
  since the $(3,3)$-coefficient divided by the $(3,1)$-coefficient
  of any matrix product of the form
  $\tilde n^{-1}\tilde m^{-1}\eta_m r$, $\tilde m\in M(F)$, $\tilde n\in N(F)$,
  is in $\OF_L$.

 \item Let $r=\begin{bmatrix}1\\w&1\\&&1&-w\\&&&1\end{bmatrix}
  \begin{bmatrix}1&\\&1\\y\varpi&z\varpi&1\\z\varpi&&&1\end{bmatrix}s_1s_2$
  with $w\in\OF/\p^n$ and $y,z\in\OF/\p^{n-1}$.

  Then $\eta_mr\notin M(F)N(F)K^\#(\P^n)$,
  since the $(3,3)$-coefficient of any product of the form
  $\tilde n^{-1}\tilde m^{-1}\eta_m r$, $\tilde m\in M(F)$, $\tilde n\in N(F)$,
  is in $\P$.

 \item Let $r = \begin{bmatrix}1&w\varpi\\&1\\&&1\\&&-w\varpi&1\end{bmatrix}s_2s_1$
  with $w\in\OF/\p^{n-1}$.

  Then $\eta_mr\notin M(F)N(F)K^\#(\P^n)$,
  since the $(4,1)$-coefficient of any product of the form
  $\tilde n^{-1}\tilde m^{-1}\eta_m r$, $\tilde m\in M(F)$, $\tilde n\in N(F)$,
  is in $\OF_L^\times$.

 \item Let $r=\begin{bmatrix}1\\w&1\\&&1&-w\\&&&1\end{bmatrix}s_1s_2s_1$ with
  $w\in\OF/\p^n$. Suppose $\eta_m r = \tilde{m}\tilde{n}k$
  with $\tilde{m} \in M(F), \tilde{n} \in N(F)$ and $k \in K^\#(\P^n)$.
  Let $A = (\tilde{m}\tilde{n})^{-1}\eta_m r$. Looking
  at the $(3,2)$ and $(3,3)$ coefficients of $A$ we get
  $\varpi^m \alpha + w \in \P^{n}$.
  If $m < n$, then we get $\alpha + w/\varpi^m \in \P$ which
  contradicts Lemma 3.1.1 (ii) of \cite{P-S}. Hence $m \geq n$, which
  implies that $w \in \P^n$. We may therefore assume that $w = 0$.
  To summarize, necessary conditions for $A \in K^\#(\P^n)$ are $m \geq n$
  and $w=0$. The following matrix identity shows that these are also
  sufficient conditions.
  \begin{equation}\label{type-6-relevance-eqn}
   \eta_m r = \begin{bmatrix}1&&&\\&&&1\\&&1&\\&-1&&\end{bmatrix}
   \begin{bmatrix}1&&&\\&1&&\\&\varpi^m \bar{\alpha}&1&\\ \varpi^m \alpha&&&1\end{bmatrix}
   \in M(F) N(F) K^\#(\P^n).
  \end{equation}

 \item Let $r = \begin{bmatrix}1&w\varpi\\&1\\&&1\\&&-w\varpi&1\end{bmatrix}
  s_2s_1s_2$ with $w\in \OF/\p^{n-1}$.

  Then $\eta_mr\notin M(F)N(F)K^\#(\P^n)$,
  since the $(3,3)$-coefficient of any product of the form
  $\tilde n^{-1}\tilde m^{-1}\eta_m r$, $\tilde m\in M(F)$, $\tilde n\in N(F)$,
  is zero.
 \item Let $r=\begin{bmatrix}1\\w&1\\&&1&-w\\&&&1\end{bmatrix}s_1s_2s_1s_2$ with
  $w\in\OF/\p^n$.

  Then $\eta_mr\notin M(F)N(F)K^\#(\P^n)$,
  since the $(3,3)$-coefficient of any product of the form
  $\tilde n^{-1}\tilde m^{-1}\eta_m r$, $\tilde m\in M(F)$, $\tilde n\in N(F)$,
  is zero.
\end{enumerate}
Let us summarize the double cosets that can possibly make a non-trivial
contribution to the integral (\ref{local-zeta-integral-2}).
\begin{align}
 &\bigcup\limits_{\substack{w \in \OF/\p^{n-1}\\z \in (\p^{n-m-1}\cap\OF)/\p^{n-1}}} K_{l,m}
  \begin{bmatrix}1&w\varpi\\&1\\&&1\\&&-w\varpi&1\end{bmatrix}
  \begin{bmatrix}1\\&1\\&&1\\&z\varpi&&1\end{bmatrix} K^\#(\p^n)
  \qquad\text{for }l,m\geq0, \label{Wsuppcoset1U}\\
 &\bigcup\limits_{\substack{w \in \OF/\p^n\\
   \varpi^m\alpha + w \in \OF_L^{\times}\\y \in (\p^{n-m-1}\cap\OF)/\p^{n-1}}}
   K_{l,m}\begin{bmatrix}1\\w&1\\&&1&-w\\&&&1\end{bmatrix}
   \begin{bmatrix}1\\&1\\yw\varpi&y\varpi&1\\y\varpi&&&1\end{bmatrix}
   s_1 K^\#(\p^n)\qquad\text{for }l,m\geq0, \label{Wsuppcoset2U}\\[1.5ex]
 &K_{l,m}s_1s_2s_1K^\#(\p^n)\qquad\text{for }l\geq0,\:m\geq n.\label{Wsuppcoset3U}
\end{align}
\subsection{Step 3: Disjointness of double cosets}
We will now investigate the overlap between double cosets in (\ref{Wsuppcoset1U}), (\ref{Wsuppcoset2U}) and (\ref{Wsuppcoset3U}). First we will consider the case $m=0$.
\subsubsection*{Equivalences among double cosets from (\ref{Wsuppcoset1U}) with $m=0$}
For $w \in \OF/\p^{n-1}$, set $\beta = c + b(\varpi w) + a(\varpi w)^2 \in \OF^{\times}$. Let $g = \begin{bmatrix}x+yb/2&yc\\-ya&x-yb/2\end{bmatrix}$ with $y = \varpi w$ and $x = c + y b/2$. Then we have the matrix identity
$$h(l,0)^{-1}\mat{g}{}{}{\det(g)\,^tg^{-1}} h(l,0) = \begin{bmatrix}1&w\varpi\\&1\\&&1\\&&-w\varpi&1\end{bmatrix} \begin{bmatrix}\beta &&&\\-a\varpi w&c&&\\&&c&a\varpi w\\&&&\beta\end{bmatrix}.$$
The rightmost matrix above is in $K^\#(\p^n)$, so that
\begin{equation}\label{Wsuppcoset1m0disjU}
 \bigcup\limits_{w \in \OF/\p^{n-1}} K_{l,0}
 \begin{bmatrix}1&w\varpi\\&1\\&&1\\&&-w\varpi&1\end{bmatrix}  K^\#(\p^n) =
 K_{l,0} K^\#(\p^n)\qquad\text{for all }l\geq0.
\end{equation}
\subsubsection*{Equivalences among double cosets from (\ref{Wsuppcoset1m0disjU})
and (\ref{Wsuppcoset2U}) with $m=0$}
Let $w \in \OF/\p^n$ be such that $\alpha+w \in \OF_L^{\times}$. Set $\beta = a + b w + c w^2$. Let $g = \begin{bmatrix}x+yb/2&yc\\-ya&x-yb/2\end{bmatrix}$ with $y=1$ and $x = -(cw + b/2)$. Then we have the matrix identity
$$h(l,0)^{-1} \mat{g}{}{}{\det(g)\,^tg^{-1}} h(l,0) \begin{bmatrix}1\\w&1\\&&1&-w\\&&&1\end{bmatrix} s_1 = \begin{bmatrix}c&&&\\-(b+cw)&-\beta&&\\&&\beta&-(b+cw)\\&&&-c\end{bmatrix}
$$
The matrix on the right hand side above is in $K^\#(\p^n)$ if $\beta \in \OF^{\times}$. We will now show that the condition $\alpha+w \in \OF_L^{\times}$ forces $\beta \in \OF^{\times}$. First observe the identity
$$
 a+bw+cw^2=-c(\alpha+w)(\alpha-(w+bc^{-1})).
$$
If $\beta \in\p$, then it would follow that
$\alpha-(w+bc^{-1})\in\p\OF_L=\P$. By Lemma 3.1.1 (ii) of
\cite{P-S}, this is impossible. It follows that indeed
$\beta\in\OF^\times$, so that \emph{all} double cosets in
(\ref{Wsuppcoset2U}) with $m=0$ are equivalent to the double coset in
(\ref{Wsuppcoset1m0disjU}).
\subsubsection*{Equivalence among double cosets from (\ref{Wsuppcoset1U}) or
(\ref{Wsuppcoset2U}) and (\ref{Wsuppcoset3U}) with $m > 0$}
Let $h_1$ be a double coset representative obtained in either (\ref{Wsuppcoset1U}) or (\ref{Wsuppcoset2U}) and let $h_2$ be a double coset representative from (\ref{Wsuppcoset3U}). Then, in either case, the double cosets are not equivalent, since, for any $r \in R(F)$ the $(2,2)$ coordinate of the matrix $h_2^{-1} h(l,m)^{-1} r h(l,m) h_1$ is in $\p$.
\subsubsection*{Equivalence among double cosets from (\ref{Wsuppcoset1U}) and
(\ref{Wsuppcoset2U}) with $m > 0$}
For $m>0$ the condition $\varpi^m\alpha + w \in \OF_L^{\times}$ in
(\ref{Wsuppcoset2U}) is equivalent to $w\in\OF^\times$. Hence let
$w \in \OF^{\times}$ and $z \in (\p^{n-m-1}\cap\OF)/\p^{n-1}$.
Let $\beta_1 = a \varpi^{2m} + b\varpi^m + c$,
$\beta_2 = a\varpi^{2m} + b\varpi^m + c w$ and $ \beta_3 = a \varpi^{2m} + b w
\varpi^m + c w^2$. We have $\beta_1, \beta_2, \beta_3 \in
\OF^\times$. Let $g =\begin{bmatrix}x+yb/2&yc\\-ya&x-yb/2\end{bmatrix}$ with
$y=\varpi^m(1-w)/\beta_3$ and $x =\beta_2/\beta_3-by/2$. Then we
have the matrix identity
\begin{eqnarray*}
&& h(l,m)^{-1} \mat{g}{}{}{\det(g)\,^tg^{-1}} h(l,m) \begin{bmatrix}1\\1&1\\&&1&-1\\&&&1\end{bmatrix}
  \begin{bmatrix}1\\&1\\z\varpi/w&z\varpi/w&1\\z\varpi/w&&&1\end{bmatrix}s_1 \\
  && \qquad = \begin{bmatrix}1\\w&1\\&&1&-w\\&&&1\end{bmatrix}
  \begin{bmatrix}1\\&1\\zw\varpi&z\varpi&1\\z\varpi&&&1\end{bmatrix}s_1 \kappa,
\end{eqnarray*}
where
$$
 \kappa = \begin{bmatrix}1&0&0&0\\\frac{c(1-w)}{\beta_3}&\frac{\beta_1}{\beta_3}&0&0\\
  \frac{cz\varpi(w^2-1)}{w\beta_3}&-\frac{\varpi^{m+1}z(w-1)(b+a\varpi^m)}
  {w\beta_3}&\frac{\beta_1}{\beta_3} & \frac{c(w-1)}{\beta_3}\\
   -\frac{\varpi^{m+1}z(w-1)(bw+a\varpi^m)}{w\beta_3}
   & -\frac{\varpi^{m+1}z(w-1)(bw+a\varpi^m(1+w))}{w\beta_3}&0&1\end{bmatrix}\in K^\#(\p^n).
$$
Hence
\begin{align}\label{Wsuppcoset1Umg0eq1}
 &\bigcup\limits_{\substack{w\in\OF/\p^n\\w\in\OF^\times\\
  z\in(\p^{n-m-1}\cap\OF)/\p^{n-1}}}
  K_{l,m}\begin{bmatrix}1\\w&1\\&&1&-w\\&&&1\end{bmatrix}
  \begin{bmatrix}1\\&1\\zw\varpi&z\varpi&1\\z\varpi&&&1\end{bmatrix}
  s_1 K^\#(\p^n) \nonumber\\
 & \qquad =\bigcup\limits_{z \in (\p^{n-m-1}\cap\OF)/\p^{n-1}} K_{l,m}
  \begin{bmatrix}1\\1&1\\&&1&-1\\&&&1\end{bmatrix}
  \begin{bmatrix}1\\&1\\z\varpi&z\varpi&1\\z\varpi&&&1\end{bmatrix}s_1 K^\#(\p^n).
\end{align}
Now let $w \in \OF/\p^{n-1}$ and $z \in (\p^{n-m-1} \cap \OF)/\p^{n-1}$. Set $\beta = c + (\varpi^{m+1} w)b + (\varpi^{m+1} w)^2 a \in \OF^{\times}$. Let $g_1 = \begin{bmatrix}x_1+y_1b/2&y_1c\\-y_1a&x_1-y_1b/2\end{bmatrix}$ with $y_1=\varpi^{m+1}w/\beta$ and $ x_1 = 1-by_1/2-a\varpi^{m+1} w y_1$.
Then we have the matrix identity
$$
 h(l,m)^{-1} \mat{g_1}{}{}{\det(g_1)\,^t(g_1)^{-1}} h(l,m)
  \begin{bmatrix}1\\&1\\&&1\\&z\varpi&&1\end{bmatrix}
 =\begin{bmatrix}1&w\varpi\\&1\\&&1\\&&-w\varpi&1\end{bmatrix}
  \begin{bmatrix}1\\&1\\&&1\\&z\varpi&&1\end{bmatrix} \kappa_1,
$$
where
$$
\kappa_1 = \begin{bmatrix}1&&&\\-a\varpi^{2m+1}w/\beta&c/\beta&&\\&a\varpi^{2+2m}wz/\beta&c/\beta&a\varpi^{2m+1}w/\beta \\a\varpi^{2+2m}wz/\beta& \varpi^{2+m}w(b+a\varpi^{m+1}w)z/\beta&&1\end{bmatrix} \in K^\#(\p^n).$$
Hence
\begin{align}\label{Wsuppcoset1Umg0eq2}
 &\bigcup\limits_{\substack{w \in \OF/\p^{n-1}\\z \in (\p^{n-m-1}\cap\OF)/\p^{n-1}}} K_{l,m}
  \begin{bmatrix}1&w\varpi\\&1\\&&1\\&&-w\varpi&1\end{bmatrix}
  \begin{bmatrix}1\\&1\\&&1\\&z\varpi&&1\end{bmatrix} K^\#(\p^n) \nonumber\\
 &\qquad=\bigcup\limits_{z \in (\p^{n-m-1}\cap\OF)/\p^{n-1}} K_{l,m}
  \begin{bmatrix}1\\&1\\&&1\\&z\varpi&&1\end{bmatrix} K^\#(\p^n).
\end{align}
We will now show that the double cosets in (\ref{Wsuppcoset1Umg0eq1}) are
all equivalent to double cosets in (\ref{Wsuppcoset1Umg0eq2}).
Given $z \in (\p^{n-m-1} \cap \OF)/\p^{n-1}$, let
$g_2 = \begin{bmatrix}x_2+y_2b/2&y_2c\\-y_2a&x_2-y_2b/2\end{bmatrix}$ with
$y_2=\varpi^m$ and $x_2 = -(c+by_2/2)$.
Then we have the matrix identity
$$
 h(l,m)^{-1} \mat{g_1}{}{}{\det(g_1)\,^t(g_1)^{-1}} h(l,m)
  \begin{bmatrix}1\\1&1\\&&1&-1\\&&&1\end{bmatrix}
  \begin{bmatrix}1\\&1\\z\varpi&z\varpi&1\\z\varpi&&&1\end{bmatrix}s_1
  = \begin{bmatrix}1\\&1\\&&1\\&z\varpi&&1\end{bmatrix} \kappa_2,
$$
where
$$\kappa_2 = \begin{bmatrix}c&&&\\-c-b\varpi^m& -c-b\varpi^m-a\varpi^{2m}&&\\ -\varpi(c+b\varpi^m)z& a\varpi^{1+2m}z&c+b\varpi^m+a\varpi^{2m} & -c-b\varpi^m\\ b\varpi^{m+1}z& \varpi^{m+1}(b+a\varpi^m)z& 0&-c\end{bmatrix} \in K^\#(\p^n).$$

We conclude that, for $m > 0$ and any $l\geq0$, the double cosets in
(\ref{Wsuppcoset1U}) and (\ref{Wsuppcoset2U}) are all contained in the union
\begin{equation}\label{Wsuppcoset1Umg0remaining}
 \bigcup\limits_{z \in (\p^{n-m-1}\cap \OF)/\p^{n-1}} K_{l,m}
 \begin{bmatrix}1\\&1\\&&1\\&z\varpi&&1\end{bmatrix} K^\#(\p^n).
\end{equation}
\subsubsection*{Equivalence among double cosets from (\ref{Wsuppcoset1Umg0remaining})
with $m > 0$}
Finally, we have to determine any equivalences amongst the double cosets
in (\ref{Wsuppcoset1Umg0remaining}). Fix $l\geq0$ and $m>0$, and let
$$
 h_1 = \begin{bmatrix}1\\&1\\&&1\\&z_1\varpi&&1\end{bmatrix}, \qquad h_2 =
 \begin{bmatrix}1\\&1\\&&1\\&z_2\varpi&&1\end{bmatrix}
$$
with $z_1,z_2 \in (\p^{n-m-1}\cap \OF)/\p^{n-1}$. We want to see if we can find $r=\mat{g}{gX}{}{\det(g)\,^tg^{-1}}\in R(F)$ such that
$$
 A = h_1^{-1}h(l,m)^{-1} r h(l,m) h_2 \in K^\#(\p^n);
$$
here, $g = \begin{bmatrix}x+yb/2&yc\\-ya&x-yb/2\end{bmatrix} \in T(F)$ and
$X = \mat{e}{f}{f}{g}$ with $e,f,g\in F$. Suppose such an $r$ exists.
Looking at the $(1,3)$, $(1,4)$, $(2,3)$ and $(2,4)$ coefficient of $A$ we get
$$
 \begin{bmatrix}x+yb/2&yc\\-ya&x-yb/2\end{bmatrix}
 \mat{e}{f}{f}{g}  \in \mat{\p^{l+2m}}{\p^{l+m}}{\p^{l+m}}{\p^l}.
$$
Looking at the  $(1,1)$, $(1,2)$, $(1,4)$ and $(3,3)$ coefficient of $A$, we see that
$$
 x \pm by/2 \in \OF^{\times}, y \in \p^m\qquad\text{and hence}\qquad
 \mat{e}{f}{f}{g}\in\mat{\p^{l+2m}}{\p^{l+m}}{\p^{l+m}}{\p^l}.
$$
Looking at the $(4,2)$ coefficient of $A$, we get
\begin{equation}\label{coeff42eq}
 (x-by/2)z_1 + \varpi^{1-l}(g(x-by/2)-afy)z_1z_2 - (x+by/2)z_2 \in \p^{n-1}.
\end{equation}
From this it follows that $\nu(z_1)=\nu(z_2)$. Using $y \in\p^m$,
it further follows that
\begin{equation}\label{4-2-coeff}
 (z_1-z_2) + \varpi^{-l}g(\varpi z_1 z_2) \in \p^{n-1}.
\end{equation}
(first add $byz_2$ to both sides of (\ref{coeff42eq}), then divide by the unit $x-by/2$).
Let $\nu(z_1) = \nu(z_2) = j$. Write $z_i = \varpi^j u_i$ for $i=1,2$, where $u_i \in \OF^{\times}$.  If $2j+1 \geq n-1$, then (\ref{4-2-coeff}) implies that
$z_1 = z_2$ which gives us that $h_1$ and $h_2$ define disjoint double cosets. If $2j+1 < n-1$, then (\ref{4-2-coeff}) implies that $u_1 - u_2 \in \p^{j+1}$.
This is a necessary condition for the coincidence of double cosets.

We will now show that it is sufficient. So, suppose that $u_1 -
u_2 \in \p^{j+1}$.  Set $g = \varpi^l(z_2-z_1)/(\varpi z_1 z_2)
\in \p^l$ and $e=f=0$. Then there is a matrix identity
$$
 h(l,m)^{-1} \mat{I_2}{X}{}{I_2} h(l,m)
 \begin{bmatrix}1\\&1\\&&1\\&z_2\varpi&&1\end{bmatrix}
 = \begin{bmatrix}1\\&1\\&&1\\&z_1\varpi&&1\end{bmatrix}
 \begin{bmatrix}1&&&\\&\frac{u_2}{u_1}&&
 \frac{z_2-z_1}{\varpi z_1z_2}\\&&1&\\&&&\frac{u_1}{u_2}\end{bmatrix},
$$
where the rightmost matrix lies in $K^\#(\p^n)$. We therefore
get the disjoint union
\begin{align*}
 \bigcup\limits_{z \in (\p^{n-m-1}\cap\OF)/\p^{n-1}}&
  K_{l,m}\begin{bmatrix}1\\&1\\&&1\\&z\varpi&&1\end{bmatrix}  K^\#(\p^n) \\
 & = \bigsqcup\limits_{z \in (\p^{n-m-1}\cap\OF\cap\p^{\big[\frac{n-1}2\big]})/\p^{n-1}}
  K_{l,m}\begin{bmatrix}1\\&1\\&&1\\&z\varpi&&1\end{bmatrix}  K^\#(\p^n) \\
 &\qquad\;\;\sqcup\bigsqcup\limits_{j=\max(n-m-1,0)}^{\big[\frac{n-3}2\big]}\;
  \bigsqcup\limits_{u\in \OF^{\times}/(1+\p^{j+1})}
  K_{l,m}\begin{bmatrix}1\\&1\\&&1\\&u\varpi^{j+1}&&1\end{bmatrix}K^\#(\p^n).
\end{align*}
\subsubsection*{Summary}
The following proposition summarizes our results in this section.
\begin{proposition}\label{disjoint-relevant}
 Let $l,m\geq0$. The following are the disjoint double cosets in
 $K_{l,m} \backslash K^H / K^\#(\p^n)$ that can possibly make a non-trivial
 contribution to the integral (\ref{local-zeta-integral-2}).
 \begin{alignat}{2}
  \bigsqcup\limits_{z \in (\p^{n-m-1}
   \cap\OF\cap\p^{\big[\frac{n-1}2\big]})/\p^{n-1}}
   &K_{l,m}\begin{bmatrix}1\\&1\\&&1\\&z\varpi&&1\end{bmatrix}  K^\#(\p^n)
   &&\text{for }l,m\geq0, \label{disjoint-cosets-support0} \\[1ex]
  \bigsqcup\limits_{j=\max(n-m-1,0)}^{\big[\frac{n-3}2\big]}\;
   \bigsqcup\limits_{u\in \OF^{\times}/(1+\p^{j+1})}
   &K_{l,m}\begin{bmatrix}1\\&1\\&&1\\&u\varpi^{j+1}&&1\end{bmatrix}  K^\#(\p^n)
   &\qquad&\text{for }l,m\geq0,\label{disjoint-cosets-support1} \\[1ex]
  &K_{l,m}s_1 s_2 s_1
   K^\#(\p^n)&&\text{for }l\geq0,\;m\geq n.\label{disjoint-cosets-support2}
 \end{alignat}
 For $n=1$ this reduces to
 \begin{alignat}{2}
  &K_{l,m}K^\#(\p)&&\text{for }l,m\geq0,\label{n=1-disjoint-support1}\\
  &K_{l,m}s_1 s_2 s_1 K^\#(\p)&\qquad
   &\text{for }l\geq0,\;m\geq1.\label{n=1-disjoint-support2}
 \end{alignat}
\end{proposition}
\section{Volume computations}\label{volume-sec}
With a view towards the integral (\ref{local-zeta-integral-2}), we will now
compute the volumes of the sets $K_{l,m} \backslash K_{l,m} A K^\#(\p^n)$,
where $A$ is one of the representatives of the disjoint double cosets in
(\ref{disjoint-cosets-support0}), (\ref{disjoint-cosets-support1}) or
(\ref{disjoint-cosets-support2}). As in Sect.\ 3.8 of \cite{P-S}, we have
\begin{equation}\label{generalvolumecalclemmaeq3new}
 \int\limits_{K_{l,m}\backslash K_{l,m}AK^\#(\p^n)}dh
 ={\rm vol}(K^\#(\p^n))
 \Big(\int\limits_{K_{l,m}\cap\big(AK^\#(\p^n)A^{-1}\big)}dt\Big)^{-1}.
\end{equation}
Note that
\begin{equation}\label{Ksharppvolumeeqnew}
 {\rm vol}(K^\#(\p^n))=\frac{q-1}{q^{3(n-1)}(q+1)(q^4-1)}
\end{equation}
from (\ref{kgeq2KsharpKHdecompeq2U})
and the fact that ${\rm vol}(K^H)=1$. Hence we are reduced to computing
\begin{equation}\label{generalvolumecalclemmaeq4new}
 V(l,m,A):=\int\limits_{K_{l,m}\cap\big(AK^\#(\p^n)A^{-1}\big)}dt.
\end{equation}
\subsection{Volume of double cosets (\ref{disjoint-cosets-support0})
and (\ref{disjoint-cosets-support1})}
In this case $A =
\begin{bmatrix}1\\&1\\&&1\\&z\varpi&&1\end{bmatrix}$ for $z \in
(\p^{n-m-1} \cap \OF)/ \p^{n-1}$. We want to find the volume of
the set $h(l,m)^{-1} R(F) h(l,m) \cap A K^\#(\p^n) A^{-1}$. Let
$\nu(z) = j$ with $j\leq n-1$.
Conjugation of $h(l,m)^{-1} R(F) h(l,m) \cap A
K^\#(\p^n) A^{-1}$ with an element of the form ${\rm diag}(1,1,u,u)$,
where $u\in\OF^\times$, leaves $R(F)$ and
$K^\#(\p^n)$ unchanged, and results in replacing $z$ by $uz$ without any
change in the volume. We may therefore assume that $z=\varpi^j$.
Since $j\leq n-1$, it is clear that
\begin{equation}\label{j-simplification}
 A K^\#(\p^n) A^{-1} \subset K^\#(\p^{j+1}).
\end{equation}
If we write an element of $R(F)$ as $tn$ with $t = \mat{x+by/2}{yc}{-ya}{x-by/2} \in T(F)$
and $n = \mat{1_2}{X}{}{1_2}$, $X = \mat{e}{f}{f}{g}$, then (\ref{j-simplification})
gives the following necessary condition for $h(l,m)^{-1} tn h(l,m) \in A K^\#(\p^n) A^{-1}$,
\begin{equation}\label{necessary condition}
 \mat{x+by/2}{yc\varpi^{-m}}{-ya\varpi^m}{x-by/2} \in
 \mat{\OF^{\times}}{\p^{j+1}}{\OF}{\OF^{\times}} \subset \GL_2(\OF)\quad\mbox{ and }\quad
 X \in \mat{\p^{2m+l}}{\p^{m+l}}{\p^{m+l}}{\p^{l}}.
\end{equation}
Set $B = A^{-1}h(l,m)^{-1} tn h(l,m)A$. We want to find further necessary conditions for $B \in K^\#(\p^n)$. Looking at the $(4,2)$ coefficient of $B$, we get
\begin{equation}\label{g-nece-condition}
 \varpi^{-l}g(x+by/2)\varpi^{2+2j} \in \p^n,
 \qquad\text{and hence}\qquad g \in \p^{n-2-2j+l}.
\end{equation}
Using the $(4,3)$ coefficient of $B$, we get
\begin{equation}\label{y-f-nece-condition}
 \varpi^l c y + \varpi^{j+1} f (x \pm by/2) \in \p^{n+m+l}.
\end{equation}
A direct computation shows that the conditions (\ref{necessary condition}), (\ref{g-nece-condition}) and (\ref{y-f-nece-condition}) are also sufficient to conclude that $B \in K^\#(\p^n)$.
Note that $\varpi^l c y + \varpi^{j+1} f (x + by/2) \in \p^{n+m+l}$ and $y \in \p^{m+j+1}$
implies that $f \in \p^{m+l}$ and $\varpi^l c y + \varpi^{j+1} f (x - by/2) \in \p^{n+m+l}$.
To summarize, the following are the necessary and sufficient conditions on $t$ and $n$ for
$h(l,m)^{-1} tn h(l,m) \in A K^\#(\p^n) A^{-1}$.
\begin{eqnarray}
&& y \in \p^{m+j+1}, \qquad x \pm by/2 \in \OF^{\times} \nonumber \\
&& e \in \p^{2m+l}, \qquad g \in \p^{n-2-2j+l} \cap \p^l, \qquad \varpi^l c y + \varpi^{j+1} f (x + by/2) \in \p^{n+m+l}. \label{nece-suff-condition}
\end{eqnarray}
For fixed values of $x, y$ satisfying the first two conditions, we are interested in
\begin{align*}
 &{\rm vol}(\{ (e,f,g) \in F^3 : e \in \p^{2m+l}, \, g \in \p^{n-2-2j+l} \cap \p^l, \,
  \varpi^l c y + \varpi^{j+1} f (x + by/2) \in \p^{n+m+l} \})  \\
 &\qquad= {\rm vol}(\{ (e,f,g) \in F^3 : e \in \p^{2m+l}, \, g \in \p^{n-2-2j+l} \cap \p^l,
  \, f \in \p^{n+m+l-j-1} - \varpi^{l-j-1} c y (x+by/2)^{-1} \}) \\
 &\qquad={\rm vol}(\{ (e,f,g) \in F^3 : e \in \p^{2m+l}, \, g \in \p^{n-2-2j+l} \cap \p^l,
  \, f \in \p^{n+m+l-j-1}\}).
\end{align*}
Note that if $j \leq \Big[\frac{n-3}2\Big]$, then $n - 2 - 2j \geq 0$,
and if $j \geq \Big[\frac{n-1}2\Big]$, then $n - 2 - 2j \leq 0$. Hence, the above volume is
\begin{eqnarray}
 q^{-2n-3m-3l+3j+3} &\mbox{ if }& j \leq \Big[\frac{n-3}2\Big]; \nonumber\\
 q^{-n-3m-3l+j+1} &\mbox{ if }& j \geq \Big[\frac{n-1}2\Big]. \label{efg-volume}
\end{eqnarray}
By an argument similar to Lemma 3.8.3 of \cite{P-S}, we get
\begin{equation}\label{t-volume}
{\rm vol}(T(F) \cap \mat{\OF^{\times}}{\p^{m+j+1}}{\OF}{\OF^{\times}})^{-1} = (1-\Big(\frac L{\p}\Big)q^{-1})q^{m+j+1}.
\end{equation}
\subsection{Volume of double coset (\ref{disjoint-cosets-support2})}
In this case, we have $A = s_1s_2s_1$ and $m \geq n$. Note that
$$
 V(l,m,s_1s_2s_1) = \int\limits_{(h(l,m)^{-1} R(F) h(l,m)) \cap
 \big(s_1s_2s_1K^\#(\p^n)(s_1s_2s_1)^{-1}\big)}dt.
$$
We have
\begin{equation}\label{compact-weyl-gp}
 s_1s_2s_1 K^\#(\p^n)(s_1s_2s_1)^{-1} =K^H\cap
 \begin{bmatrix}\OF^{\times}&\OF&\OF&\p^n\\ \p^n&\OF^{\times}&\p^n&\p^n\\
  \p&\OF&\OF^{\times}&\p^n\\ \OF&\OF&\OF&\OF^{\times}\end{bmatrix}.
\end{equation}
We have to find the intersection of this compact group with $h(l,m)^{-1} R(F) h(l,m)$. Set
$$
 L_1 := \mat{\OF^{\times}}{\OF}{\p^n}{\OF^{\times}} \subset \GL_2(\OF),
  \quad N_1 := \{X \in \mat{\OF}{\p^n}{\p^n}{\p^n} :\, ^tX = X\} \subset F^3.
$$
Then $L_1$ and $N_1$ are the upper left and upper right blocks of
(\ref{compact-weyl-gp}), respectively. Write a given element of
$R(F)$ as $tn$ with $t\in T(F)$ and $n\in U(F)$. If
$n=\mat{1_2}{X}{}{1_2}$, then by arguments similar to those in
Sect.\ 3.8 of  \cite{P-S}, we see that $tn$ lies in
$s_1s_2s_1K^\#(\p^n)(s_1s_2s_1)^{-1}$ if and only if
\begin{equation}\label{volcomp1cond1aeqnew}
 \mat{\varpi^{-m}}{}{}{1}t\mat{\varpi^m}{}{}{1}\in L_1
\end{equation}
and
\begin{equation}\label{volcomp1cond2aeqnew}
\mat{\varpi^{-2m-l}}{}{}{\varpi^{-m-l}}X\mat{1}{}{}{\varpi^m}
 \in N_1.
\end{equation}
It follows that
\begin{align*}
 &{\rm vol}\big(\{X\in F^3:\:  \mat{\varpi^{-2m-l}}{}{}{\varpi^{-m-l}}X\mat{1}{}{}{\varpi^m}\,
\in N_1\}\big)\\
 &\qquad={\rm vol}\big(\{X\in F^3:\:X\in
  \mat{\varpi^{2m+l}}{}{}{\varpi^{m+l}}N_1
  \mat{1}{}{}{\varpi^{-m}}\}\big)\\
 &\qquad=q^{-3m-3l}{\rm vol}(N_1) =q^{-3m-3l-2n}.
\end{align*}
Let
\begin{equation}\label{generalvolumecalclemma2eq2new}
  T_m =\{t\in T(F):\:
  \mat{\varpi^{-m}}{}{}{1}t\mat{\varpi^m}{}{}{1}\in L_1\}
  =T(F)\cap\mat{\varpi^m}{}{}{1}L_1\mat{\varpi^{-m}}{}{}{1}.
\end{equation}
So far, we have $V(l,m,s_1s_2s_1)^{-1} = q^{3m+3l+2n} {\rm vol}(T_m)^{-1}$.
By an argument similar to Lemma 3.8.4 of \cite{P-S}, and using $m\geq n$, we get
\begin{equation}\label{t-volume-1}
  {\rm vol}(T_{m})^{-1}=\Big(1-\Big(\frac L\p\Big)q^{-1}\Big)q^m.
\end{equation}
The following proposition summarizes the volume computations in this section.
\begin{proposition}\label{volumecompsummarylemma}
 \begin{enumerate}
  \item Let $m \geq 0$. Let $A = \begin{bmatrix}1\\&1\\&&1\\&z\varpi&&1\end{bmatrix}$
   for $z \in (\p^{n-m-1} \cap \OF)/ \p^{n-1}$ and set $\nu(z) = j$.
   If $j \leq \Big[\frac{n-3}2\Big]$, then
   \begin{equation}\label{type1volumeformula1}
    V_j^{l,m} := \int\limits_{K_{l,m}\backslash K_{l,m}AK^\#(\p^n)}dh
    = \frac{q-1}{q^{3(n-1)}(q+1)(q^4-1)} (1-\Big(\frac L{\p}\Big)q^{-1})q^{2n+4m+3l-2j-2},
   \end{equation}
   and if $j \geq \Big[\frac{n-1}2\Big]$, then
   \begin{equation}\label{type1volumeformula2}
    V^{l,m} := \int\limits_{K_{l,m}\backslash K_{l,m}AK^\#(\p^n)}dh
     = \frac{q-1}{q^{3(n-1)}(q+1)(q^4-1)} (1-\Big(\frac L{\p}\Big)q^{-1})q^{n+4m+3l}.
   \end{equation}
  \item For all $m \geq n$,
   \begin{equation}\label{type6volumeformula}
    V^{l,m}_{s_1s_2s_1} := \int\limits_{K_{l,m}\backslash
    K_{l,m}s_1s_2s_1K^\#(\p^n)}dh = \frac{q-1}{q^{3(n-1)}(q+1)(q^4-1)}
     \Big(1-\Big(\frac L\p\Big)q^{-1}\Big)q^{4m+3l+2n}.
   \end{equation}
  \item In particular, for $n=1$,
   \begin{alignat*}{2}
    \int\limits_{K_{l,m}\backslash K_{l,m}K^\#(\p)}dh
    &=\frac{q-1}{(q+1)(q^4-1)} (1-\Big(\frac L{\p}\Big)q^{-1})q^{4m+3l+1}&\qquad&(m\geq0),\\
    \int\limits_{K_{l,m}\backslash K_{l,m}s_1s_2s_1K^\#(\p)}dh
    &=\frac{q-1}{(q+1)(q^4-1)} \Big(1-\Big(\frac L\p\Big)q^{-1}\Big)q^{4m+3l+2}&&(m>0).
   \end{alignat*}
 \end{enumerate}
\end{proposition}
Note that the right hand side of (\ref{type1volumeformula2}) is independent of $j$.
This will play an important role in the evaluation of the zeta integral.
\section{Main local theorem}\label{main-thm-sec}
In this section we will calculate the integral (\ref{local-zeta-integral-2}). From Proposition \ref{disjoint-relevant}, we have
\begin{align}
 Z(s)&= \sum\limits_{l,m \geq 0} B(h(l,m)) \Big(\sum\limits_{z \in (\p^{n-m-1}
  \cap\OF\cap\p^{\big[\frac{n-1}2\big]})/\p^{n-1}} W^\#(\eta h(l,m) A(z), s) V^{l,m}
  \nonumber \\
 &\hspace{20ex} + \sum\limits_{j=\max(n-m-1,0)}^{\Big[\frac{n-3}2\Big]}\;
  \sum\limits_{u \in \OF^{\times}/(1+\p^{j+1})} W^\#(\eta h(l,m) A(\varpi^{j+1}u),s)
  V_j^{l,m}\Big) \nonumber \\
 &+ \sum\limits_{l \geq 0, m \geq n} B(h(l,m)) W^\#(\eta h(l,m)s_1s_2s_1,s)
  V^{l,m}_{s_1s_2s_1} \label{zeta-integral-formula}
\end{align}
where $A(z) = \begin{bmatrix}1\\&1\\&&1\\&z\varpi&&1\end{bmatrix}$.
By (\ref{Wsharpformulaeq}), (\ref{type-1-relevance-eqn}) and (\ref{type-6-relevance-eqn})
we get
\begin{align}
 W^\#(\eta h(l,m) A(z), s) &=|\varpi^{2m+l}|^{3(s+\frac 12)} \omega_{\pi}(\varpi^{-2m-l})
  \omega_{\tau}(\varpi^{-m-l}) W^{(0)}(\mat{\varpi^l}{0}{\varpi z}{1}),
  \label{W-sharp-formula-1} \\
 W^\#(\eta h(l,m)s_1s_2s_1,s) &=|\varpi^{2m+l}|^{3(s+\frac 12)}
  \omega_{\pi}(\varpi^{-2m-l}) \omega_{\tau}(\varpi^{-m-l})
   W^{(0)}(\mat{}{\varpi^l}{-1}{}).\label{W-sharp-formula-2}
\end{align}
Set $C_{l,m} := |\varpi^{2m+l}|^{3(s+\frac 12)} \omega_{\pi}(\varpi^{-2m-l})
\omega_{\tau}(\varpi^{-m-l})$. Substituting (\ref{W-sharp-formula-1}) and
(\ref{W-sharp-formula-2}) into (\ref{zeta-integral-formula}), we get
\begin{align}\label{zeta-integral-formula-2}
 Z(s)&=\sum\limits_{l \geq 0} B(h(l,0)) C_{l,0} W^{(0)}(\mat{\varpi^l}{0}{0}{1}) V^{l,0}
 \nonumber \\
 &+ \sum\limits_{l \geq 0, m \geq 1} B(h(l,m)) C_{l,m} V^{l,m}
  \Big( \sum\limits_{z \in (\p^{n-m-1}\cap\OF\cap\p^{\big[\frac{n-1}2\big]})/\p^{n-1}}
  W^{(0)}(\mat{\varpi^l}{0}{\varpi z}{1})\Big) \nonumber \\
 &+ \sum\limits_{l \geq 0, m \geq 1} B(h(l,m)) C_{l,m}
  \Big(\sum\limits_{j=\max(n-m-1,0)}^{\Big[\frac{n-3}2\Big]}\;
  \sum\limits_{u \in \OF^{\times}/(1+\p^{j+1})} W^{(0)}(\mat{\varpi^l}{0}{\varpi^{j+1}u}{1})
  V^{l,m}_j\Big) \nonumber \\
 &+ \sum\limits_{l \geq 0, m \geq n} B(h(l,m)) C_{l,m} W^{(0)}(\mat{}{\varpi^l}{-1}{})
  V^{l,m}_{s_1s_2s_1}.
\end{align}
If $n=1$, then the inner sum over $z$ in the second term above
reduces to just $z=0$, and the third term above is not present. In
this case, the integral $Z(s)$ was computed in Theorem 3.9.1 of \cite{P-S}.

From now on we will assume that $n \geq 2$. As the following lemma shows,
the fact that the representation $\tau$ has conductor $\p^n$ implies that the
middle two expressions in formula (\ref{zeta-integral-formula-2}) are zero.
\begin{lemma}\label{W-vectors-zero}
 Let $m \geq 1$ and $n \geq 2$.
 \begin{enumerate}
  \item For any $g \in \GL_2(F)$,
   $$
    \sum\limits_{z \in (\p^{n-m-1}\cap\OF\cap\p^{\big[\frac{n-1}2\big]})/\p^{n-1}}
     W^{(0)}(g \mat{1}{0}{\varpi z}{1}) = 0.
   $$
  \item For $2j+2 < n$ and any $z$ with $\nu(z) = j$,
   $$
    W^{(0)}(\mat{\varpi^l}{0}{\varpi z}{1}) = 0.
   $$
 \end{enumerate}
\end{lemma}
{\bf Proof:} i) Let $t=\max(n-m-1, 0, \big[\frac{n-1}2\big])$.
We have $\p^{n-m-1}\cap\OF\cap\p^{\big[\frac{n-1}2\big]} = \p^t$ and, since $m \geq 1$ and $n \geq 2$, we see that $t+1 < n$. Define $\hat{W}(g) = \sum\limits_{z \in \p^{t+1}/\p^n} W^{(0)}(g \mat{1}{0}{z}{1}) \in V_{\tau}$. Then, by definition $\hat{W}$ is right invariant under the group $\mat{1}{0}{\p^{t+1}}{1}$. Since $W^{(0)}$ is right invariant under $K^{(1)}(\p^n)$, we see that $\hat{W}$ is right invariant under $\mat{\OF^\times}{}{}{\OF^\times}$. The matrix identity $\mat{1}{x}{}{1} \mat{1}{}{z}{1} = \mat{1}{}{(1+xz)^{-1}z}{1} \mat{1+xz}{x}{}{(1+xz)^{-1}}$ and $1+xz \in \OF^\times$ implies that $\hat{W}$ is also right invariant under the group $\mat{1}{\OF}{}{1}$. Since $\mat{1}{0}{\p^{t+1}}{1}, \mat{\OF^\times}{}{}{\OF^\times}$ and $\mat{1}{\OF}{}{1}$ generate $K^{(1)}(\p^{t+1})$, it follows that $\hat{W}$ is a vector in $V_{\tau}$ that is right invariant under $K^{(1)}(\p^{t+1})$. Since $\tau$ has level $\p^n$ and $t+1 < n$,
this implies $\hat{W} = 0$, as claimed.
\nll
ii) Let $z_1, z_2$ be such that $\nu(z_1) = \nu(z_2) = j$ and $z_1/z_2 \in 1+\p^{j+1}$.
Consider the matrix identity
$$
 \mat{\varpi^l}{}{}{1} \mat{1}{}{\varpi z_1}{1}
 = \mat{1}{\frac{\varpi^l(z_2-z_1)}{\varpi z_1z_2}}{}{1} \mat{\varpi^l}{}{}{1}
 \mat{1}{}{\varpi z_2}{1} \mat{\frac{z_1}{z_2}}
 {\frac{(z_2-z_1)}{\varpi z_1z_2}}{}{\frac{z_2}{z_1}}.
$$
Since the additive character $\psi$ is trivial on $\OF$ and the rightmost matrix
is in $K^{(1)}(\p^n)$, it implies that
\begin{equation}\label{W-vectors-zeroeq1}
 W^{(0)}(\mat{\varpi^l}{}{}{1} \mat{1}{}{\varpi z}{1}) = W^{(0)}(\mat{\varpi^l}{}{}{1}
 \mat{1}{}{\varpi zu}{1})
\end{equation}
for every $u \in 1+\p^{j+1}$ and $z\in\OF$ with $\nu(z)=j$
(we have essentially derived the well-definedness of the third sum in
(\ref{zeta-integral-formula-2})). Writing $u = 1+b \varpi^{j+1}$ with $b \in \OF$ and
integrating both sides of (\ref{W-vectors-zeroeq1}) with respect to $b$, we get
\begin{align*}
 W^{(0)}(\mat{\varpi^l}{}{}{1} \mat{1}{}{\varpi z}{1})
  &=\int\limits_{\OF} W^{(0)}(\mat{\varpi^l}{}{}{1} \mat{1}{}{\varpi z}{1}
   \mat{1}{}{\varpi z b \varpi^{j+1}}{1})db \\
 &=\int\limits_{\OF} W^{(0)}(\mat{\varpi^l}{}{}{1} \mat{1}{}{\varpi z}{1}
  \mat{1}{}{b \varpi^{2j+2}}{1}) db.
\end{align*}
This last expression is zero, since $2j+2 < n$ and
$\tilde{W}(g) := \int\limits_{\OF} W^{(0)}(g \mat{1}{}{b \varpi^{2j+2}}{1}) db \in V_{\tau}$
is right invariant under $K^{(1)}(\p^{2j+2})$. This concludes the proof.\qed
\nl
Using this lemma, (\ref{zeta-integral-formula-2}) now becomes
\begin{align}\label{zeta-integral-formula-3}
 Z(s) &=\sum\limits_{l\geq0}B(h(l,0))|\varpi^{l}|^{3(s+\frac 12)}\omega_{\pi}(\varpi^{-l})
  \omega_{\tau}(\varpi^{-l}) W^{(0)}(\mat{\varpi^l}{0}{0}{1}) V^{l,0} \nonumber \\
 &+\sum\limits_{l \geq 0, m \geq n} B(h(l,m))
  |\varpi^{2m+l}|^{3(s+\frac 12)} \omega_{\pi}(\varpi^{-2m-l})
  \omega_{\tau}(\varpi^{-m-l}) W^{(0)}(\mat{}{\varpi^l}{-1}{})
  V^{l,m}_{s_1s_2s_1}.
\end{align}
Since $\mat{}{1}{\varpi^n}{}$ normalizes $K^{(1)}(\p^n)$, the vector
$W'(g) := W^{(0)}(g \mat{}{1}{\varpi^n}{})$ is another element of $V_{\tau}$
that is right invariant under $K^{(1)}(\p^n)$. Since the space of vectors in $V_{\tau}$
right invariant under $K^{(1)}(\p^n)$ is one dimensional, there is a constant $c \in \C$
such that $W^{(0)} = c W'$ (one can check that $c^{-2} = \omega_{\tau}(\varpi^n)$).
Hence,
\begin{equation}\label{W-0-equation}
 W^{(0)}(\mat{}{\varpi^l}{-1}{}) = c W^{(0)}(\mat{\varpi^{l+n}}{}{}{-1}) = c
 W^{(0)}(\mat{\varpi^{l+n}}{}{}{1}).
\end{equation}
This shows that in order to evaluate
(\ref{zeta-integral-formula-3}) we need the formula for the
new-vector of $\tau$ in the Kirillov model. The possibilities for our
generic, irreducible, admissible representation $\tau$ of $\GL_2(F)$
with unramified central character and conductor $\p^n$, $n\geq2$, are as follows.
Either $\tau$ is a principal series representation $\chi_1\times\chi_2$,
where $\chi_1$ and $\chi_2$ are ramified characters of $F^\times$ (with $\chi_1\chi_2$
unramified); or $\tau=\chi\St_{\GL(2)}$, a twist of the Steinberg representation
by a ramified character $\chi$ (such that $\chi^2$ is unramified); or $\tau$ is
supercuspidal. In each case the newform in the Kirillov model is given by
the characteristic function of $\OF^\times$; see, e.g.,  \cite{Sc}.
It follows that all the terms in (\ref{W-0-equation}) are zero.
The integral (\ref{zeta-integral-formula-3}) reduces to
\begin{equation}\label{zeta-integral-formula-4}
 Z(s) = V^{0,0} = \frac{q-1}{q^{3(n-1)}(q+1)(q^4-1)} (1-\Big(\frac L{\p}\Big)q^{-1})q^n.
\end{equation}
Thus, we have proved the following result.
\begin{theorem}\label{localintegralmaintheoremn>1}
 Let $\pi$ be an unramified, irreducible, admissible representation of $\GSp_4(F)$
 (not necessarily with trivial central character), and let
 $\tau$ be an irreducible, admissible, generic representation of $\GL_2(F)$
 with unramified central character and conductor $\p^n$ with $n \geq 2$.
 Let $Z(s)$ be the integral (\ref{local-zeta-integral}), where $W^\#$ is the function
 defined in Sect.\ 2, and $B$ is the spherical Bessel
 function defined in Sect.\ 1 (ix). Then
 \begin{equation}\label{localintegralmaintheoremn>1eq}
  Z(s)=\frac{q-1}{q^{3(n-1)}(q+1)(q^4-1)} (1-\Big(\frac L{\p}\Big)q^{-1})q^n.
 \end{equation}
\end{theorem}
{\bf Remark:} For any unramified, irreducible, admissible representation $\pi$ of $\GSp_4(F)$ and any of the representations $\tau$ of $\GL_2(F)$ mentioned in the
theorem we have $L(s, \pi \times \tau) = 1$.
Hence, up to a constant, the integral $Z(s)$ represents the $L$-factor
$L(s, \pi \times \tau)$.
\section{Global integral and special value of $L$-function}\label{global-section}
Let $\Gamma_2 = \SSp_4(\Z)$. For a positive integer $l$ denote by
$S_l(\Gamma_2)$ the space of Siegel cusp forms of degree $2$ and
weight $l$ with respect to $\Gamma_2$.
Let  $\Phi \in S_l(\Gamma_2)$ be a Hecke eigenform. It has a Fourier expansion
$$
 \Phi(Z) = \sum\limits_{S > 0}a(S,\Phi)e^{2 \pi i \tr(SZ)},
$$
where $S$ runs through all symmetric, semi-integral, positive
definite matrices of size two. Let us make the following two
assumptions about the function $\Phi$.
\begin{description}
 \item[Assumption $1$:] $a(S,\Phi) \neq 0$ for some $S =
  \mat{a}{b/2}{b/2}{c}$ such that $b^2-4ac = -D < 0$ where $-D$
  is the discriminant of the imaginary quadratic field $L := \Q(\sqrt{-D})$.
 \item[Assumption $2$:] The weight $l$ is a multiple of $w(-D)$, the
  number of roots of unity in $\Q(\sqrt{-D})$. We have
  $$
   w(-D) = \left\{\begin{array}{ll}
    4 & \hbox{ if } D=4, \\
    6 & \hbox{ if } D=3, \\
    2 & \hbox{ otherwise. }\end{array}\right.
  $$
\end{description}
We lift the function $\Phi$ to a function $\phi_\Phi$ on $H(\A) = \GSp_4(\A)$,
where $\A$ is the ring of adeles of $\Q$, in a standard way; see (141) in \cite{P-S}.
Let $V_{\Phi}$ be the automorphic representation generated by $\phi_{\Phi}$, and let
$\pi_\Phi \cong \otimes'_p \pi_{p}$ be an irreducible component.
Let $\Lambda = \otimes \Lambda_p$ be a character of $L^\times$ depending on $S$ as
constructed in Sect.\ 5.1 of \cite{P-S}. 

Let $N = \prod p^{n_p} \in \N$. Denote the space of Maa{\ss} cusp
forms of weight $l_1 \in \Z$ with respect to $\Gamma_0(N)$ by
$S_{l_1}^M(N)$. A function $f \in S_{l_1}^M(N)$ has the Fourier
expansion
\begin{equation}\label{Maass form four exp}
 f(x+iy) = \sum\limits_{n \neq 0} a_n W_{{\rm sgn}(n)\frac{l_1}2,
\frac{ir}2}(4 \pi|n|y)e^{2 \pi i nx},
\end{equation}
where $W_{\nu,\mu}$ is a classical Whittaker function and
$(\Delta_{l_1}+\lambda) f=0$ with $\lambda = 1/4 + (r/2)^2$. Here
$\Delta_{l_1}$ is the Laplace operator defined in Sect.\ 5.3 of
\cite{P-S}. Let $f \in S_{l_1}^M(N)$ be a  Hecke eigenform.

If $ir/2 = (l_2-1)/2$ for some integer $l_2 > 0$, then the
cuspidal, automorphic representation of $\GL_2(\A)$ constructed
below is holomorphic at infinity of lowest weight $l_2$. In this
case we make the additional assumptions that $l_2\leq l$ and
$l_2\leq l_1$, where $l$ is the weight of the Siegel cusp form
$\Phi$. Starting from $f$, we obtain another Maa{\ss} form $f_l
\in S_{l}^M(N)$ by applying the raising and lowering operators as
in (147) of \cite{P-S}. Define a function $\hat{f}$ on $\GL_2(\A)$
by
\begin{equation}\label{maass form lift to group}
 \hat{f}(\gamma_0 m k_0) = \Big(\frac{\gamma i +\delta}{|\gamma i +\delta|}\Big)^{-l}
 f_l\Big(\frac{\alpha i + \beta}{\gamma i + \delta}\Big),
\end{equation}
where $\gamma_0 \in \GL_2(\Q)$, $m =
\mat{\alpha}{\beta}{\gamma}{\delta} \in \GL_2^{+}(\R)$, $k_0 \in
\prod\limits_{p | N}K^{(1)}(\p^{n_p}) \prod\limits_{p \nmid
N}\GL_2(\Z_p)$. Here, for $p | N$ we have $K^{(1)}(\p^{n_p}) =
\GL_2(\Q_p) \cap\mat{\Z_p^{\times}}{\Z_p}{\p^{n_p}}{\Z_p^\times}$
with $\p = p\Z_p$. Let $\tau_f \cong \otimes'_p \tau_{p}$ be the irreducible,
cuspidal, automorphic representation of $\GL_2(\A)$ generated by
$\hat f$. As in Sect.\ 5.2 of \cite{P-S}, define an Eisenstein
series on $\GU(2,2;L)(\A)$ by
\begin{equation}\label{eisenstein series definition}
 E_{\Lambda}(g,s) = \sum\limits_{\gamma \in P(\Q) \backslash
 G(\Q)}f_{\Lambda}(\gamma g,s)
\end{equation}
where $f_\Lambda$ is as defined in (154) of \cite{P-S} from
$\hat{f}$. We consider the global integral
\begin{equation}\label{global integral calculation}
 Z(s,\Lambda) = \int\limits_{Z_H(\A)H(\Q)\backslash H(\A)}E_{\Lambda}(h,s)\bar{\phi}(h)dh.
\end{equation}
Now, applying Theorem 3.7 from \cite{Fu}, Theorem 4.4.1 from \cite{P-S} and Theorem \ref{localintegralmaintheoremn>1}, we get
\begin{thm}\label{main-global-thm}
 Let $\Phi \in S_l(\Gamma_2)$ be a cuspidal Siegel eigenform of degree $2$ and even
 weight $l$ satisfying the two assumptions above.
 Let $L=\Q(\sqrt{-D})$, where $D$ is as in Assumption 1.
 Let $N = \prod p^{n_p}$ be a  positive integer. Let $f$ be a Maa{\ss}
 Hecke eigenform of weight $l_1 \in \Z$ with respect to $\Gamma_0(N)$. If $f$ lies in a
 holomorphic discrete series with lowest weight $l_2$, then assume that $l_2 \leq l$.
 Then the integral (\ref{global integral calculation}) is given by
 \begin{equation}\label{integral-l-fn-formula}
  Z(s,\Lambda) = \kappa_{\infty} \kappa_N \frac{L(3s+\frac 12, \pi_{\Phi} \times
  \tau_f)}{\zeta(6s+1) L(3s+1, \tau_f\times \AI(\Lambda))},
 \end{equation}
 where
 \begin{align*}
  \kappa_{\infty} &=\frac12\overline{a(\Lambda)}c(1)\pi D^{-3s-\frac l2}\,
   (4 \pi)^{-3s+\frac 32-l}\,\frac{\Gamma(3s+l-1+\frac{ir}2) \Gamma(3s+l-1-\frac{ir}2)}
   {\Gamma(3s+\frac{l+1}2)}, \\
   \kappa_N &=\prod\limits_{p | N} \frac{p-1}{p^{3(n_p-1)}(p+1)(p^4-1)} (1-\Big(\frac L{p}\Big)p^{-1})p^{n_p}   (1-p^{-6s-1})^{-1} \prod\limits_{p^2 | N} L_p(3s+1, \tau_p \times \AI(\Lambda_p)).
 \end{align*}
 Here, the non-zero constant $c(1)$ is given by Eqn. (148) of \cite{P-S},
 the non-zero constant $a(\Lambda)$ is defined in
 Sect.\ 5.1 of \cite{P-S}, and
 $$
  \Big(\frac Lp\Big)=\left\{\begin{array}{l@{\qquad\text{if }p}l}
  -1&\text{ is inert in }L,\\
  0&\text{ ramifies in }L,\\
  1&\text{ splits in }L.
  \end{array}\right.
 $$
 The quantity $\frac{ir}2$ is as in (\ref{Maass form four exp}).
\end{thm}

One would like to know if the local $L$-function $L_p(s, \tau_p \times \AI(\Lambda_p))$ is
$1$. If $p$ is an odd prime, then the only case where the $L$-function
$L_p(s, \tau_p \times \AI(\Lambda_p))$ is not identically $1$ is when $p | D$,
$\nu_p(N) = 2$ and $\tau_p$ is a certain induced representation or a certain
twist of the Steinberg representation.
The main difficulty is that if $\nu_p(N) = 2$, then it is not clear if the corresponding
representation $\tau_p$ is induced or special or supercuspidal.

Let $\Gamma^{(2)}(N) := \{g \in \SSp_4(\Z) : g \equiv
1 \pmod{N} \}$ be the principal congruence subgroup of $\SSp_4(\Z)$. Let us denote the space of all
Siegel cusp forms of weight $l$ with respect to
$\Gamma^{(2)}(N)$ by $S_l(\Gamma^{(2)}(N))$.
For a Hecke eigenform $\Phi \in S_l(\Gamma^{(2)}(N))$, let $\Q(\Phi)$ be the subfield of $\C$ generated by all the Hecke eigenvalues of $\Phi$. Let $S_l(\Gamma^{(2)}(N), \Q(\Phi))$ be the subspace of $S_l(\Gamma^{(2)}(N))$ consisting of cusp forms whose Fourier coefficients lie in $\Q(\Phi)$. For more details on this space we refer to Sect.\ 5.4 of \cite{P-S}. Note that all the arguments of Sect.\ 5.4 of \cite{P-S} are valid, with minor modifications, if we remove the restriction that $N$ is square-free and use the definition (\ref{U-compact defn}) for $K^\#(\P^n)$. Hence, we get the following special value result.

\begin{thm}\label{special values thm}
 Let $\Phi$ be a cuspidal Siegel eigenform of weight $l$ with respect to
 $\Gamma_2$ satisfying the two assumptions above, let $D$ be as in Assumption 1 and let 
 $\Phi \in S_l(\Gamma^{(2)}(N), \Q(\Phi))$. Let $N = \prod p^{n_p}$ be odd such that
 if $p | D$ then $n_p \neq 2$. Let $\Psi$ be a normalized, holomorphic,
 cuspidal eigenform of weight $l$ with respect to $\Gamma_0(N)$. Then
 \begin{equation}\label{special-values-eqn}
  \frac{L(\frac l2 - 1, \pi_{\Phi} \times \tau_{\Psi})}{\pi^{5l-8}
  \langle \Phi,\Phi \rangle \langle \Psi,\Psi \rangle_1} \in \bar{\Q},
 \end{equation}
 where $\langle \Phi,\Phi \rangle$ and $\langle \Psi,\Psi \rangle_1$ denote
 the Petersson inner products of $\Phi$ and $\Psi$, respectively, and $\bar{\Q}$
 denotes the algebraic closure of $\Q$ in $\C$.
\end{thm}
Note that we have the above restriction on $N$ because if $p | D$ and $n_p = 2$, then we do not know if the term $L_p((l-1)/2, \tau_p \times \AI(\Lambda_p))$ in $\kappa_N$ from Theorem \ref{main-global-thm} is algebraic or not.

\end{document}